\documentclass[12pt,a4paper]{amsart}
\usepackage{amsmath,amsfonts,amsthm,amssymb,latexsym,mathabx}
\usepackage{xspace}
\usepackage{xcolor}
\usepackage{enumitem}
\usepackage{graphicx}
\usepackage{tabularx}
\usepackage{longtable}
\usepackage{epsfig}
\usepackage{hhline}
\usepackage{float}
\usepackage{hyperref}

\newcommand{\soc}{\mathrm{soc}}

\newcommand{\Aut}{\mathrm{Aut}}
\newcommand{\Lg}{\mathrm{L}}
\newcommand{\ol}{\overline}
\newcommand{\Ug}{\mathrm{U}}

\newcommand{\PSL}{\mathrm{PSL}}

\newcommand{\PSp}{\mathrm{PSp}}

\theoremstyle{plain}
\newtheorem{theorem}[]{Theorem}[section] 
\newtheorem{lemma}[theorem]{Lemma} 
\newtheorem{proposition}[theorem]{Proposition}
\newtheorem{corollary}[theorem]{Corollary}

\theoremstyle{plain}
\newtheorem{mtheorem}[]{Theorem}
\newtheorem{mcorollary}[mtheorem]{Corollary}

\theoremstyle{definition}
\newtheorem*{definition}{Definition}

\theoremstyle{plain}
\newtheorem*{conjecture}{Conjecture}

\date{\today}
\title{On finite groups with soluble centralisers}
\author[V. Grazian]{Valentina Grazian}
\address{Valentina Grazian: Department of Mathematics \\ University of Padova \\ Via Trieste 63, 35121 Padova (PD) \\ Italy}
\email{valentina.grazian@unipd.it}

\author[C. Monetta]{Carmine Monetta}
\address{Carmine Monetta: Department of Mathematics \\ University of Salerno \\
via Giovanni Paolo II 132, 84084 Fisciano (SA)\\ Italy}
\email{cmonetta@unisa.it}

\author[G. Tracey]{Gareth Tracey}
\address{Gareth Tracey: Mathematics Institute \\ University of Warwick \\ Coventry CV4 7AL \\ United Kingdom}
\email{gareth.tracey@warwick.ac.uk}


\begin{document}

\begin{abstract}
We classify finite groups in which the centralisers of certain non-central elements are soluble. This includes a full structural description of groups whose non-central element centralisers are all soluble, and a reduction theorem for the case in which all non-central $\pi$-elements have soluble centralisers, for a suitable collection $\pi$ of primes. Our results yield further descriptions under mild local conditions and have applications to groups with soluble involution centralisers, as well as to questions concerning non-commuting graphs.
\end{abstract}

\maketitle

\section{Introduction}\label{sec:intro}

The classification of finite groups according to the structure of their centralisers is a long-standing and challenging problem in group theory. In particular, the task of describing all finite groups in which every centraliser is soluble remains open and continues to attract attention from different perspectives. Over time, several related problems have been investigated as simplified versions of this question. For example, groups in which all centralisers are metabelian (that is, soluble of derived length two) or nilpotent have been considered (see \cite{CAS15, SUZ61}).  

Another line of research focuses on the centralisers of specific types of elements. For instance, Problem 2.74 from the Kourovka Notebook \cite{KM} asks for a description of finite groups in which all non-central involutions have soluble centralisers. At present, the existing results essentially concern the case where the Sylow 2-subgroups are abelian (cf. \cite{GOR2}). This limitation highlights how intricate the general problem is and provides further motivation to study groups in which only certain non-central elements, possibly distinguished by their order, have soluble centralisers.  

In particular, even partial information about the sections that such groups may admit represents meaningful progress toward understanding their overall structure. In this paper, we contribute to this direction by investigating finite groups $G$ in which the centralisers $C_G(x)$ of certain non-central elements $x$ of $G$ are soluble.

First of all we tackle the case in which the centralisers of \emph{all} non-central elements of $G$ are soluble: in this instance we characterize the quotient $G/R(G)$, where $R(G)$ denotes the soluble radical of $G$. In what follows, $\pi(n)$ denotes the set of prime divisors of a positive integer $n$; and we write $\pi(G)$ for $\pi(|G|)$, when $G$ is a finite group.

\begin{mtheorem}\label{thm:2ndMain}
Let $G$ be a finite group with the property that $C_G(x)$ is soluble for every non-central element $x$ in $G$. Set $\mathcal{Q}:=\{2^{p_1},3^{p_2},p\mid p_i,p\text{ odd primes, }p\equiv 0,\pm 2\pmod{5}\}$. Then one of the following holds.
\begin{enumerate}
    \item[\upshape(i)] $G$ is soluble.
    \item[\upshape(ii)] $G/R(G)$ is almost simple with socle one of the following:
    \begin{enumerate}
        \item[\upshape(1)] $\mathrm{M}_{11}$, $\mathrm{M}_{22}$, $A_5$, $A_6$, $A_7$, $\mathrm{G}_2(3)$, $\PSp_4(3)$, ${}^2\mathrm{F}_4(2)'$;
    \item[\upshape(2)] $\Lg^{\pm}_4(q)$, $\PSp_4(q)$, $\mathrm{G}_2(q)$ or ${}^2\mathrm{F}_4(q)$ with $q\in\mathcal{Q}$;
    \item[\upshape(3)] $\Lg_2(q)$, $\Lg^{\pm}_3(q)$,  ${}^2\mathrm{B}_2(q)$, ${}^2\mathrm{G}_2(q)$.   
\end{enumerate} 
    \item[\upshape(iii)] For one of the sets of primes $\pi$ in Table \ref{tab:Main}, $G/R(G)$ is a $\pi$-group and all non-abelian composition factors of $G$ lie in $S(\pi)$.
\end{enumerate}

\renewcommand{\arraystretch}{1.5}
\begin{table}[H]
    \centering
    \begin{tabular}{c|c}
      $\pi$    &  $S(\pi)$\\
      \hline
      $\pi({}^2\mathrm{B}_2(2^p))\cup\{p\}$, $p$ an odd prime & $ {}^2\mathrm{B}_2(2^p)$\\
      
      $\bigcup_{q\in\mathcal{Q}_0}\pi(\Aut(\Lg_2(q)))$ with $\mathcal{Q}_0\subseteq \mathcal{Q}$ of size at most $3$ & $\Lg_2(q), q\in\mathcal{Q}_0$\\
      $\{2,3,5\}$ & $ A_5, A_6, \PSp_4(3)$\\
      $\{2,3,7\}$ & $ \Lg_3(2),  \Lg_2(8),  \Ug_3(3)$\\
      $\{2,3,13\}$ & $\Lg_3(3)$ \\\hline
      \end{tabular}  \bigskip
    \caption{Primes and composition factors from Theorem \ref{thm:2ndMain}(iii).}
    \label{tab:Main}
\end{table}
\end{mtheorem}

It is plausible that part (iii) of Theorem \ref{thm:2ndMain} could be eliminated, but doing so appears to require theoretical advances, potentially beyond cohomology, rather than further refinements of existing techniques. In this sense, (iii) reflects a natural endpoint of the available theory.

\begin{conjecture}\label{conj:Main}
Let $G$ be a finite insoluble group with the property that $C_G(x)$ is soluble for every non-central element $x$ in $G$, and let $R(G)$ be the soluble radical of $G$. Then $G/R(G)$ is almost simple with socle one of $\mathrm{M}_{11}$, $\mathrm{M}_{22}$, $A_n$ ($n\le 7$), $\Lg_2(q)$, $\Lg_3(3)$, $\Lg_3(4)$, $\mathrm{PSp}_4(3)$, $\Ug_3(3)$, $\Ug_4(2)$, $\Ug_5(2)$, $\mathrm{G}_2(3)$, ${}^2\mathrm{B}_2(q)$ or ${}^2\mathrm{F}_4(2)'$.   
\end{conjecture}

We next extend our analysis to the case where the solubility condition on centralisers is imposed only on the $\pi$-elements of the group, for a given set of primes~$\pi$. Our first aim is to characterize the non-abelian simple $\pi$-sections of $G$.

To do this, we have to fix some notation.
\begin{definition}
Let $\pi$ be a set of primes and let $p$ be the smallest prime in $\pi$. For a prime power $q$, set $\delta(q):=2$ if $q$ is odd, and $\delta(q):=1$ if $q$ is even. We then define the following collections of finite simple groups.

\renewcommand{\arraystretch}{1.5}
\begin{table}[H]
    \centering
    \begin{tabular}{c|c}
    \hline
$X_{\mathrm{Alt}}(\pi)$ & $\{A_5,A_6,A_7\}$ if $3\in\pi$\\
                         & $\{A_5,A_6,A_7,A_8\}$ if $2\in\pi,3\not\in\pi$\\
                         & $\{A_n\text{ : }n\le p+4\}$ if $2\not\in\pi,3\not\in\pi$\\ \hline
$X_{\mathrm{Class}}(\pi)$                          & $\{\Lg^{\pm}_n(q),\Omega^{\epsilon}_n(q),\PSp_{\delta(q)n}(q) \text{ : }n\le 7,q\text{ a prime power}\}$ if $3\in\pi$\\
                          & $\{\Lg^{\pm}_n(q),\Omega^{\epsilon}_n(q),\PSp_{\delta(q)n}(q) \text{ : }n\le 8,q\text{ a prime power}\}$ if $2\in\pi$, $3\not\in\pi$\\ 
 & $\{\Lg^{\pm}_n(q),\Omega^{\epsilon}_n(q),\PSp_{\delta(q)n}(q) \text{ : }n\le p+4,q\text{ a prime power}\}$ if $2,3\not\in\pi$\\ \hline
$X_{\mathrm{Exc}}$ & All exceptional groups of Lie type\\ \hline
$X_{\mathrm{Spor}}(\pi)$ & $\{ S \text{ a sporadic group : } \text{$S$ and $\overline{\pi}$ are as in Table \ref{tab:spor} and }\overline{\pi} \text{ contains }\pi\}$\\
\hline
\end{tabular} \bigskip
    \caption{Sets of simple groups associated to a set of primes $\pi$, with $p$ the smallest prime in $\pi$.}
    \label{tab:Xpi}
\end{table}

\begin{table}[H]
    \begin{tabular}{|c|c|}
    \hline
      $S$   &  $\overline{\pi}$\\
      \hline
      $\mathrm{M}_{11}$   &  $2,3,5,11$\\
      $\mathrm{M}_{12}$   &  $3,5,11$\\
      $\mathrm{M}_{22}$   &  $2,3,5,7,11$\\
      $\mathrm{M}_{23}$   &  $5,7,11,23$\\
      $\mathrm{M}_{24}$   &  $5,7,11,23$\\
      $\mathrm{HS}$   &  $7,11$\\
        $\mathrm{J}_{2}$   &  $7$\\
      $\mathrm{Co}_{1}$   &  $11,13,23$\\
      $\mathrm{Co}_{2}$   &  $7,11,23$\\
      \hline
    \end{tabular}
    \quad
        \begin{tabular}{|c|c|}
        \hline
      $S$   &  $\overline{\pi}$\\
      \hline
      $\mathrm{Co}_{3}$   &  $7,11,23$\\
      $\mathrm{McL}$   &  $5,7,11$\\
      $\mathrm{Suz}$   &  $7,11,13$\\
      $\mathrm{He}$   &  $17$\\
      $\mathrm{HN}$   &  $11,19$\\
      $\mathrm{Th}$   &  $5,13,19,31$\\
      $\mathrm{Fi}_{22}$   &  $7,11,13$\\
      $\mathrm{Fi}_{23}$   &  $11,13,17,23$\\
      $\mathrm{Fi}_{24}'$   &  $11,13,17,23,29$\\
      \hline
    \end{tabular}
    \quad
        \begin{tabular}{|c|c|}
        \hline
      $S$   &  $\overline{\pi}$\\
      \hline
      $\mathrm{B}$   &  $11,13,17,19,23,31,47$\\
      $\mathrm{M}$   &  $19,23,29,31,41,47,59,71$\\
      $\mathrm{J}_1$   &  $3,5,7,11,19$\\
      $\mathrm{O'N}$   &  $5,7,11,19,31$\\
      $\mathrm{J}_3$   &  $5,17,19$\\
      $\mathrm{Ru}$   &  $7,13,29$\\
      $\mathrm{J}_{4}$   &  $11,23,29,31,37,43$\\
      $\mathrm{Ly}$   &  $7,11,31,37,67$\\
      \hline
    \end{tabular}
     \bigskip
    \caption{Sporadic groups in the collection $X_{\mathrm{Spor}}(\pi)$.}
    \label{tab:spor}
\end{table}
\end{definition}

\begin{mtheorem}\label{thm:Main}
Let $G$ be a finite group, and let $\pi$ be a set of primes. Suppose that $C_G(x)$ is soluble for every non-central $\pi$-element $x$ in $G$. Then the non-abelian sections of $G$ lie in $$X_{\mathrm{Alt}}(\pi)\cup X_{\mathrm{Class}}(\pi)\cup X_{\mathrm{Exc}} \cup X_{\mathrm{Spor}}(\pi). $$
\end{mtheorem}

A slightly sharper formulation of Theorem \ref{thm:Main} is possible.
Indeed, by invoking Proposition \ref{prop:PSp} one may further restrict the symplectic composition factors that can occur. We have chosen the present formulation in order to emphasise the uniform structure of the conclusion and to avoid introducing additional case-specific notation.

 We point out that the proof of Theorem \ref{thm:Main} relies on cohomological arguments concerning groups that arise as non-split extensions of its Frattini subgroup by a simple group, showing the diversity of techniques required to tackle the problem.

Theorems \ref{thm:2ndMain} and \ref{thm:Main} show that if $G$ is a finite group in which \emph{all} centralisers of non-central $\pi$-elements are soluble, then the non-abelian simple sections of $G$ are either $\pi'$-groups; or they have bounded order in the alternating cases, and bounded twisted rank in the Lie type cases. In particular, the degree of the largest alternating section of $G$ (sometimes referred to as the \emph{thickness} of $G$) is bounded in terms of $\pi$. 
\begin{mcorollary}\label{cor:degree}
Let $G$ be a finite group, and let $\pi$ be a set of primes. Suppose that $C_G(x)$ is soluble for every non-central $\pi$-element $x$ in $G$, and let $d$ be the degree of the largest alternating section of $G$. Then $d\le p+4+2\delta_{2,p}$ for all $p\in\pi$.  
\end{mcorollary}

The bound of the thickness of $G$ in Corollary \ref{cor:degree} is of course best possible, since alternating groups of degree $p+5+2\delta_{2,p}$ contain an element of order $p$ centralised by a subgroup isomorphic to $A_5$.

We also remark here that in Corollaries \ref{cor:simple2} and \ref{cor:simple}, we give the exact lists of finite simple groups with soluble centralisers of all non-trivial involutions and all non-trivial elements, respectively.

\medskip
Our methods allow us to derive a stronger conclusion to Theorem \ref{thm:Main} if the group in question satisfies a mild local condition $(\ast)_\pi$:

\begin{definition}
Let $\pi$ be a set of primes. We will say that a finite group $G$ has property $(\ast)_{\pi}$ if $G$ is either $\pi$-soluble or for each prime $r\in \pi$, one of the following holds for a Sylow $r$-subgroup $R$ of $G$:
\begin{enumerate}
    \item[\upshape(a)] $R$ is $r$-central, that is, all elements of $R$ of order $r$ are contained in the center $Z(R)$; or
    \item[\upshape(b)] $G$ has a non-abelian simple section $S$ with a Sylow $r$-subgroup of the same nilpotency class as $R$.
\end{enumerate}
\end{definition}

\begin{mtheorem}\label{thm:MainClass}
Let $G$ be a finite group with property $(\ast)_{\pi}$, for a set of primes $\pi$. Suppose that $C_G(x)$ is soluble for every non-central $\pi$-element $x$ in $G$. Then either $G$ is $\pi$-soluble, or $G/Z(G)=D.A$, where $D$ is a soluble $(\pi\cap \pi(A))'$-group and $A$ is an almost simple group in which all non-trivial $\pi$-elements have soluble centraliser. In particular, $\soc(A)$ is as described in Theorem \ref{thm:Main}.
\end{mtheorem}

    

The case $\pi=\{2\}$ is of particular interest. The following result answers the question 2.74 in \cite{KM} for groups having $2$-central Sylow $2$-subgroups, extending the work of Gorenstein \cite{GOR2}.

 \begin{mcorollary}\label{cor:MainClass}
Let $G$ be a finite group with the property that $C_G(x)$ is soluble for every non-central involution $x$ in $G$, and the Sylow $2$-subgroups of $G$ are $2$-central. Then either $G$ is soluble, or $G/Z(G)$ has shape $G/Z(G)=D.A$, where $D$ has odd order and $A$ is almost simple with all involution centralisers soluble. In particular, $\mathrm{soc}(A)$ is one of $\mathrm{M}_{11}$, $\mathrm{M}_{22}$, $A_n$ ($n\le 8$), $\Lg_2(q),\Lg_3(2^f),\Ug_3(2^f),\Lg_3(3)$ $\mathrm{PSp}_4(3),\Ug_3(3),\Ug_4(2), \Ug_5(2),\mathrm{G}_2(3)$, ${}^2\mathrm{B}_2(q)$, ${}^2\mathrm{F}_4(2)'$.  \end{mcorollary}

We conclude with a reduction theorem concerning a question of Abdollahi, Akbari and Maimani \cite{AAM}. This question belongs to the graph-theoretic approach to group theory, which investigates finite groups through the structure of graphs built from their algebraic relations. For a group $G$, we denote by $\Gamma_G$ the \emph{non-commuting graph} of $G$, that is, the simple and undirected graph whose vertex set is $G\setminus Z(G)$, and with an edge between $x$ and $y$ if and only if they do not commute. In \cite[Question 1.2]{AAM}, the authors ask for group theoretic properties $\mathcal{P}$ that can be distinguished by the non-commuting graph. More precisely, for such a property $\mathcal{P}$, and a $\mathcal{P}$-group $G$, it is asked whether or not any group $H$ with $\Gamma_H\cong \Gamma_G$ needs to have property $\mathcal{P}$: for more details on problems of this type see for instance \cite{GLM} and \cite{GM}. Our result gives the structure of a minimal counterexample in the case where the property in question is solubility.   
\begin{mcorollary}\label{cor:NCgraph}
Suppose that there exist finite groups $G$ and $H$ with $G$ insoluble, $H$ soluble, and $\Gamma_G\cong \Gamma_H$. Let $G$ be a group of minimal order with this property. Then $|G|=|H|$ and each of the following holds.
\begin{enumerate}
    \item[\upshape(i)] $G$ is as in Theorem \ref{thm:2ndMain}(ii) or (iii); and
    \item[\upshape(ii)] $H$ contains elements $a$ and $b$ such that $C_H(a)\cap C_H(b)=Z(H)$.
\end{enumerate}
\end{mcorollary}
From a graph-theoretic viewpoint, condition (ii) above can be interpreted as asserting that $\Gamma_H$ has domination number at most $2$. We remark that the proof of (ii) does not require the Classification of Finite Simple Groups. 

\vspace{2mm}
\subsection*{Outline of the approach and methods.}
Our proofs combine cohomological, structural, and computational techniques. A central theme is the analysis of groups of the form \(G=\Phi(G).S\) with \(S\) non-abelian simple. This naturally leads to the study of Frattini extensions and derivations from \(S\) into irreducible \(\mathbb{F}_p[S]\)-modules. Vanishing results for the low-degree cohomology groups \(\mathrm{H}^1(S,V)\) and \(\mathrm{H}^2(S,V)\) are used to establish precise lifting criteria for elements and their centralisers (cf. Lemma \ref{lem:derivation} and Corollary \ref{cor:strategy}). These results make it possible to transfer the existence of insoluble centralisers from the simple quotient \(S\) to the full group \(G\), even when \(G\) has a nontrivial Frattini kernel. 
The argument then proceeds by induction on the order of \(G\) and by reduction to almost simple quotients, through a detailed \(\pi\)-local analysis. Another key ingredient is the use of embeddings of alternating groups into classical and sporadic simple groups, which provides effective bounds on natural module dimensions and links the alternating and Lie-type cases. In addition, explicit computations with the algebra system \textsc{Magma} were used to verify cohomology dimensions, automorphism structures, and centraliser correspondences in exceptional or small-degree cases, complementing the theoretical reductions. 
Overall, Sections~3--5 develop a unified framework that combines these cohomological lifting techniques, alternating-group embeddings, and computational checks. The resulting methods not only establish our main theorems, but also provide general tools for studying finite groups where solubility of centralisers imposes strong global restrictions on the group’s composition factors.

\bigskip
\subsection*{Layout of the paper.} Section 2 is devoted to the analysis of finite simple groups $S$ in which the centralisers of all non-trivial $\pi$-elements are soluble, with a complete characterisation when $\pi = \{2\}$ and $\pi = \pi(S)$. In Section 3, we prove Theorem \ref{thm:Main}, establishing structural properties for a broader class of groups. Section 4 deals with Theorem \ref{thm:MainClass} and contains the proofs of Corollaries \ref{cor:MainClass} and \ref{cor:NCgraph}. Finally, Section 5 is devoted to the proof of Theorem \ref{thm:2ndMain}.

\subsection*{Notation}

We follow the notation of the \emph{Atlas of Finite Groups}~\cite{Atlas1985}.

We use standard notation for subgroup operators: \(O_p(G)\) denotes the largest normal \(p\)-subgroup of \(G\); \(O_{p'}(G)\) the largest normal \(p'\)-subgroup; \(\soc(G)\) the socle of \(G\); \(R_\pi(G)\) the \(\pi\)-soluble radical; \(F(G)\) and \(F^*(G)\) the Fitting and generalized Fitting subgroups, respectively; and \(\Phi(G)\) the Frattini subgroup.

For extensions, \(H.K\) denotes a (possibly non-split) extension of \(H\) by \(K\), while \(H \circ K\) denotes a central product of \(H\) and \(K\).

If \(\mathcal{Q}\) is a set of positive integers, then \(\pi(\mathcal{Q})\) denotes the set of all prime divisors of the elements of \(\mathcal{Q}\); if \(\mathcal{Q}=\{q\}\), we write \(\pi(q)\) for \(\pi(\{q\})\). As usual, if \(H\) is a finite group, we write \(\pi(H)=\pi(|H|)\).

\section{The simple cases}
The property of having all centralisers of non-trivial $\pi$-elements soluble is preserved by subgroups. With this in mind, we start with a lemma concerning the embedding of alternating groups in finite simple groups.

\begin{lemma}\label{lem:AltEmbeddings}
Let $S$ be a finite simple group of Lie type, with natural module $V$ of dimension $n\geq 5$ if $S$ is classical. Then $S$ contains a subgroup isomorphic to the alternating group $A_m$, where
\begin{enumerate}
    \item[\upshape(i)] $m=n+1$ if $S\in\{\Lg_n(q),\Omega^{\circ}_n(q)\}$ or $S=\Ug_n(q)$ with $n$ odd.
    \item[\upshape(ii)] $m=n$ if $n$ is even and either $S\in\{\Ug_n(q),\Omega^{\pm}_n(q)\}$, or $S=\PSp_n(q)$ with $q$ even.
    \item[\upshape(iii)] $m=n/2$ if $S=\PSp_n(q)$ with $q$ odd.
    \item[\upshape(iv)] $m$ is as in Table \ref{tab:EXC} if $S$ is exceptional.
\end{enumerate}  
\end{lemma}

\begin{table}[H]
    \begin{tabular}{c|c|c|c|c|c|c|c|c|c}
    \hline
       ${}^2\mathrm{B}_{2}(q)$ & ${}^2\mathrm{G}_{2}(q)$ & ${}^3\mathrm{D}_{4}(q)$ & $\mathrm{G}_{2}(q)$ & ${}^2\mathrm{F}_{4}(q)$ &  $\mathrm{F}_{4}(q)$   & 
        ${}^2\mathrm{E}_6(q)$   &  
      $\mathrm{E}_6(q)$   &  
      $\mathrm{E}_7(q)$   & 
      $\mathrm{E}_8(q)$\\
      $2$ &  $3$ &  $5$ &  $5$ &  $5$  &  $8$ &  $10$ &  $10$   &  $10$  &  $10$\\ 
      \hline
       \end{tabular}\bigskip
    \caption{Degrees of some alternating subgroups of exceptional groups.}
    \label{tab:EXC}
\end{table}

\begin{proof}
Via its action on the deleted permutation module of dimension $n$, the alternating group $A_{n+1}$ embeds as a subgroup of $\Lg_n(q)$. If $n$ is odd, then this action in fact yields an embedding $A_{n+1}\hookrightarrow \Omega^{\circ}_{n}(q)$. It follows that $A_{n+1}$ is a subgroup of $\Omega^{\circ}_n(q)\le \Lg^{\pm}_n(q)$ in this case, and $A_{n+1}<\PSp_{n+1}(q)$ if $q$ is even. 

Now suppose that $n$ is even, and fix $\epsilon\in\{\pm\}$. 
Note that $\Omega_{n-1}(q)$ embeds as a subgroup of $\Omega_n^{\epsilon}(q)$; that $\Omega_n^{\epsilon}(q)$ embeds in $\Ug_n(q)$; and that $\Omega_{n-1}(q)\cong \PSp_{n-2}(q)$ if $q$ is even.
Hence 
$A_n$ is a subgroup of $\Ug_n(q)$ and $\Omega^{\pm}_n(q)$ for all $q$. If $q$ is odd and $S=\PSp_n(q)$, then $S$ contains a copy of its Weyl group $W(S)\cong S_2\wr S_{n/2}$, whence a copy of $S_{n/2}$.

Finally, we can deduce (iv) by inspection of the known maximal subgroups of the finite exceptional groups $S$ of Lie type. These are available in \cite{Craven, Craven2, KleidmanRee,Kleidman,Malle, suzmax} for $S\neq \mathrm{E}_8(q)$. For $S=\mathrm{E}_8(q)$, we just use the embedding $\mathrm{E}_7(q)\le \mathrm{E}_8(q)$ and \cite{Craven2}.
\end{proof}

The following result concerns alternating sections in sporadic simple groups. It will be used later in the paper and can verified by inspecting of the Atlas of Finite Groups \cite{ATLAS}.
\begin{lemma}\label{lem:sporalt}
 Let $S$ be a sporadic simple group. Then the maximal degree $n$ of an alternating section of $S$ is given in Table \ref{tab:sporalt}.   
\end{lemma}

\begin{table}[H]
    \begin{tabular}{c|c|c|c|c|c|c|c|c|c|c|c|c}
    \hline
       $\mathrm{M}_{11}$ & $\mathrm{M}_{12}$ & $\mathrm{M}_{22}$ & $\mathrm{M}_{23}$ & $\mathrm{M}_{24}$ &  $\mathrm{HS}$   & 
        $\mathrm{J}_{2}$   &  
      $\mathrm{Co}_{1}$   &  
      $\mathrm{Co}_{2}$   & 
      $\mathrm{Co}_{3}$   &  
      $\mathrm{McL}$   &  
      $\mathrm{Suz}$   &  
      $\mathrm{He}$\\
      $6$ &  $6$ &  $7$ &  $8$ &  $8$  &  $8$ &  $6$ &  $9$   &  $8$  &  $8$ &  $8$ &  $7$ &  $7$\\ 
      \hline
       \end{tabular}
       \quad
      \begin{tabular}{c|c|c|c|c|c|c|c|c|c|c|c|c}
    \hline  
      $\mathrm{HN}$   &  
      $\mathrm{Th}$   &  
      $\mathrm{Fi}_{22}$   &  
      $\mathrm{Fi}_{23}$ &
      $\mathrm{Fi}_{24}'$   &  
    $\mathrm{B}$   &  
      $\mathrm{M}$   &  
      $\mathrm{J}_1$   & 
      $\mathrm{O'N}$   &  
      $\mathrm{J}_3$   &  
      $\mathrm{Ru}$   & 
      $\mathrm{J}_{4}$   &  
      $\mathrm{Ly}$\\ 
    $12$  &  $9$  &  $10$  &  $10$  &  $12$ &  $12$ &  $12$ &  $5$ &  $7$   &  $6$  &  $8$ &  $8$ &  $11$\\
      \hline
      \end{tabular}\bigskip
    \caption{Maximal degrees of alternating sections of sporadic groups}
    \label{tab:sporalt}
\end{table}

We can now determine the finite simple groups in which the centralisers of all non-trivial $\pi$-elements are soluble.
\begin{lemma}\label{lem:simple}
Let $S$ be a non-abelian finite simple group, and let $\pi$ be a set of primes dividing $|S|$. Suppose that $C_S(x)$ is soluble for every non-trivial $\pi$-element $x$ in $S$. Define 
$X'_{\text{Class}}(\pi):=\{\Lg^{\pm}_n(q),\Omega^{\epsilon}_n(q),\PSp_{\delta(n,q)n}(q) \text{ : }n\le p+4,q\text{ a prime power}\}$.
Then 
$$S\in X_{\mathrm{Alt}}(\pi)\cup X_{\mathrm{Class}}'(\pi)\cup X_{\mathrm{Exc}} \cup X_{\mathrm{Spor}}(\pi).$$
Furthermore, if $S$ is of Lie type in characteristic $\ell$ and $\ell\in\pi$, then $S$ has twisted rank $1$ or $S\in\{\Lg_3(q),\Ug_3(q),\Lg_4(2),\Lg_4(3),\Ug_4(2),\Ug_4(3),\PSp_4(2)',\PSp_4(3),\Ug_5(2),\mathrm{G}_2(3),{}^2\mathrm{F}_4(2)'\}$. 
\end{lemma}
\begin{proof}
If $S$ is a sporadic group, then the result can be checked by inspection of the Atlas of finite groups \cite{ATLAS}. If $S$ is an alternating group of degree at least $m:=p+5+2\delta_{2,p}$, then $S$ contains a copy of $A_m$, and the centraliser of a $p$-cycle in $A_m$ is insoluble if $p$ is odd, while the centraliser of a product of transpositions in $A_m=A_9$ is insoluble in the case $p=2$. 

So we may assume that $S=X_r(q)$ is a finite group of Lie type of twisted rank $r$. Let $\ell$ be the defining characteristic of $S$, so that $q$ is a power of $\ell$. Also, write $S=X_{\sigma}$, where $X$ is a simple algebraic group and $\sigma:X\rightarrow X$ is a Steinberg endomorphism. Write $\gamma$ for the graph automorphism involved $\gamma$ in $\sigma$, so that $d:=|\gamma|\in\{1,2,3\}$.
Also, let $\Pi$ be a set of fundamental roots for $X$, and for $\alpha\in \Pi$, write $U_{\hat{\alpha}}$ for the associated twisted root subgroup, as defined in \cite[Section 2.3]{GLS3}. Since $\gamma$ is induced by a graph automorphism $\hat{\gamma}$ of the Dynkin diagram of $X$, we may define $\Pi(\alpha)$ to be the $\hat{\gamma}$-orbit containing $\alpha$ in $\Pi$. 

Let $\alpha,\beta\in \Pi$ such that $\Pi(\alpha)\neq \Pi(\beta)$, and for each $\alpha'\in \Pi(\alpha)$ and $\beta'\in\Pi(\beta)$, $\alpha'$ and $\beta'$ are non-adjacent in the Dynkin diagram for $X$. Then the Chevalley relations for $S$ \cite[Theorem 2.4.5]{GLS3} show that for then $U_{\hat{\alpha}}$ commutes with $\langle U_{\pm\hat{\beta}}\rangle$. Since $C_S(u)$ is soluble for all $u\in U_{\hat{\alpha}}$, we deduce that $r\le 3$. By inspection of the finite groups of Lie type of twisted rank at most $3$, we see that the centralisers of all $\ell$-elements of $S$ are soluble if and only if $S\in\{\Lg_2(q),\Lg_3(q),\Ug_3(q),\Lg_4(2),\Lg_4(3),\Ug_4(2),\Ug_4(3),\PSp_4(3),\Ug_5(2),\mathrm{G}_2(3),{}^2\mathrm{B}_2(q)$, ${}^2\mathrm{F}_4(2)',$  ${}^2\mathrm{G}_2(q)\}$. 

Suppose next that $s\in \pi$ with $s\neq\ell$. If $s=2$, then by inspection of \cite[Table 4.5.1]{GLS3} we see that $S\in\{\Lg_2(q),\Lg_3(3),\Ug_3(3),\PSp_4(3),\mathrm{G}_2(3)\}$.
So assume that $s>2$. The result then follows from Lemma \ref{lem:AltEmbeddings} and the first paragraph above. 
\end{proof}

Note that Lemma \ref{lem:simple} gives the second part of Theorem \ref{thm:Main}. 
We conclude this section with the complete description of the non-abelian finite simple groups $S$  having soluble centralisers of all non-trivial $\pi$-elements, first for $\pi = \{2\}$ and then for $\pi = \pi(S)$.

\begin{corollary}\label{cor:simple2}
Let $S$ be a non-abelian simple group, and suppose that the centraliser of every non-trivial involution in $S$ is soluble. Then $S$ is one of $\mathrm{M}_{11}$, $\mathrm{M}_{22}$, $A_n$ ($n\le 8$), $\Lg_2(q),\Lg_3(2^f),\Ug_3(2^f),\Lg_3(3)$ $\mathrm{PSp}_4(3),\Ug_3(3),\Ug_4(2), \Ug_5(2),\mathrm{G}_2(3),{}^2\mathrm{B}_2(q),{}^2\mathrm{F}_4(2)'$. Furthermore, each of these groups has the property that the centraliser of every non-trivial involution is soluble.        
\end{corollary}
\begin{proof}
If $S$ is alternating or sporadic, then the result follows immediately from Lemma \ref{lem:simple}. So assume that $S$ has Lie type. Applying Lemma \ref{lem:simple} and its proof, we see that the only possibilities are that $S\in\{\Lg_2(q),\Lg_3(2^f),\Ug_3(2^f),\Lg_4(2),\Ug_4(2),\Ug_5(2),{}^2\mathrm{B}_2(q),{}^2\mathrm{F}_4(2)'$, $\Lg_3(3),$ $\Ug_3(3),\PSp_4(3),\mathrm{G}_2(3)\}$.
\end{proof}

\begin{corollary}\label{cor:simple}
Let $S$ be a non-abelian simple group, and suppose that the centraliser of every non-trivial element of $S$ is soluble. Then $S$ is one of $A_n$ ($n\le 7$), $\Lg_2(q),\Lg_3(3),\Lg_3(4),$ $\mathrm{PSp}_4(3),\Ug_3(3),\Ug_4(2), \Ug_5(2),\mathrm{G}_2(3),{}^2\mathrm{B}_2(q),\mathrm{M}_{11}$, $\mathrm{M}_{22},{}^2\mathrm{F}_4(2)'$. Furthermore, each of these groups has the property that the centraliser of every non-trivial element is soluble.        
\end{corollary}
\begin{proof}
By inspection of the groups from Corollary \ref{cor:simple2}, we see that $S$ has all of its proper centralisers soluble if and only if $S$ is one of the listed cases. 
\end{proof}

\section{Lifting centralisers and the proof of Theorem \ref{thm:Main}}

The key obstacle in the proof of our main theorems concerns the idea of lifting centralisers in finite groups $G$. Namely, if we know properties of the centralisers of $G/N$ for $N\unlhd G$, can we deduce properties of the centralisers in $G$? 
Our first two preparatory lemmas address this question for certain normal subgroups $N$.


\begin{lemma}\label{lem:derivation}
Let $G$ be a finite group, let $E\unlhd G$, write $S=G/E$, and write $\rho:G\rightarrow S$ for the natural projection. Assume there exists $1\neq x\in S$ such that $(|x|,|E|)=1$ and $C_S(x)$ is insoluble. Then there exists $\hat{x}\in \rho^{-1}(x)\subseteq G$ with $|\hat{x}|=|x|$ and $C_G(\hat{x})$ projects onto $C_S(x)$ modulo $E$. In particular, $C_G(\hat{x})$ is insoluble.
\end{lemma}
\begin{proof}
Since $(|x|,|E|)=1$, the Schur-Zassenhaus theorem implies that there exists an element $\hat{x}$ of $\pi^{-1}(x)$ of order $|x|$. Then $(|\hat{x}|,|E|)=1$, so $C_G(\hat{x})E/E=C_{G/E}(E\hat{x})=C_S(x)$.
The result follows.     
\end{proof}

\begin{lemma}\label{lem:Z1}
Let $G$ be a finite group, $Z$ be a subgroup of $Z(G)$, and $p$ be prime. Suppose that $G/Z$ has a non-central element of $p$-power order with an insoluble centraliser. Then $G$ has a non-central element of $p$-power order with an insoluble centraliser. 
\end{lemma}
\begin{proof}
Write $r=|Z(G)|$, and let $Zy$ be a non-central element of $G/Z$ of $p$-power order with an insoluble centraliser $C/Z$. We may assume that $y\in G$ has order $p^a$. Suppose that $x$ is an element of $C$ of order coprime to $p$. Then $[x,y]\in Z$, so $[x,y,y]=[x,y,x]=1$. Thus, $[x,y]=1$ since $x$ and $y$ have coprime orders. Writing $|C|_p$ for the $p$-part of $|C|$, we have therefore shown that all elements $c$ of $C$ have the property that $c^{|C|_p}\in C_G(y)$. Let $D:=\langle c^{|C|_p}\text{ : }c\in C\rangle$. Then $D\unlhd C$, and $D\le C_G(y)$. Since $C/D$ has $p$-power exponent and $C$ is insoluble, we deduce that $D$ is insoluble. Thus, $C_G(y)\geq D$ is insoluble.  
\end{proof}

\begin{corollary}\label{cor:Z1}
Let $G$ be a finite group of shape $G=R.E$ where $R$ is nilpotent and $E=G/R$ is a direct product of $t>1$ non-abelian simple groups. Let $p$ be a prime dividing $|E|$. Suppose that either $O_p(G)=1$, or that there is a $G$-normal series $1=V_0< \hdots< V_k=O_p(G)=O_p(R)$ such that $V_i/V_{i-1}$ is centralised by $G$ for all $i$. Then $G$ has a non-central $p$-element with an insoluble centraliser.
\end{corollary}
\begin{proof}
Note first that since $G/R$ is a direct product of at least $2$ non-abelian simple groups and $p$ divides $|G/R|$, $G/R$ has a $p$-element with an insoluble centraliser. If $O_p(G)=1$, then $R$ is a $p'$-group, so the result follows from Lemma \ref{lem:derivation}, applied with $E:=R$.

So assume that $O_p(G)>1$, and that $1=V_0<\hdots<V_k=O_p(G)$ is a normal series with $G$ acting trivially on the factors $V_i/V_{i-1}$. If $k=1$, then the result therefore follows from Lemma \ref{lem:Z1}, applied with $Z:=V_1=O_p(G)$. For $k>1$, we can apply induction and Lemma \ref{lem:Z1} to $G/V_1$.  
\end{proof}

The following consequence of the results above will be very useful in the proofs of our main theorems. We remind the reader that if $G$ is a group acting on another (not necessarily abelian) group $E$, then a \emph{derivation} is a map $\delta:G\rightarrow E$ satisfying $\delta(gh)=\delta(g)\delta(h)^{g^{-1}}$ for all $g,h\in G$. For example, if $e\in E$, then $\delta_e:G\rightarrow E, g\mapsto [e,g^{-1}]$ defines a derivation. We define $\mathrm{Der}(G,E):=\{\delta\text{ : }\delta\text{ a derivation}\}$ and $\mathrm{InnDer}(G,E):=\{\delta_e\text{ : }e\in E\}$. 

The map $\delta\rightarrow \{\delta(g)g\text{ : }g\in G\}$ gives a bijection from $\mathrm{Der}(G,E)$ to the set of complements to $E$ in $E\rtimes G$. Note also that if $E$ is abelian, then $\mathrm{Der}(G,E)$ is a group, and $\mathrm{H}^1(G,E)=\mathrm{Der}(G,E)/\mathrm{InnDer}(G,E)$.

\begin{corollary}\label{cor:strategy}
Let $G$ be a finite group, let $E$ be a soluble normal subgroup of $G$, and let $p$ be a prime. Suppose that $S:=G/E$ has a non-central element of order $p$ with an insoluble centraliser $C$, and that one of the following holds.
\begin{enumerate}
    \item[\upshape(1)] $p$ does not divide $|E|$.
    \item[\upshape(2)] $C$ has an insoluble $p'$-subgroup.
    \item[\upshape(3)] $E=\Phi(G)$, $p$ divides $|E|$ and there is a set $\mathcal{S}$ of insoluble subgroups of $S$ such that for each non-trivial irreducible $\mathbb{F}_p[S]$-module $V$, either some member of $\mathcal{S}$ centralises a non-zero element of $V$, or $\mathrm{H}^2(S,V)=\mathrm{H}^1(D,V\hspace{-0.2cm}\downarrow_D)=0$ for all $D\in\mathcal{S}$.   
\end{enumerate}
Then $G$ contains a non-central element of order $p$ with an insoluble centraliser. 
\end{corollary}
\begin{proof}
If $p$ does not divide $E$, then the result follows from Lemma \ref{lem:derivation}, so we just need to prove that if (2) or (3) holds, then $G$ contains a non-central element of order $p$ with an insoluble centraliser. Throughout, write $\rho:G\rightarrow S$ for the natural projection, and let $z$ be a non-central element of $S$ of order $p$ with $C_S(z)=C$.

Suppose first that $C$ has an insoluble $p'$-subgroup $D$. Let $D_0$ be a subgroup of $G$ minimal with the property that $D_0\rho=D$. Then $D_0\cap E\le \Phi(D_0)$, so $\pi(|D_0|)=\pi(|D|)$. In particular, $D_0$ is a $p'$-group. We will prove, by induction on $|E|$, that $D_0$ centralises an element $x$ of $G$ of $p$-power order with $x\rho=z$. This will of course suffice to prove the corollary in this case. By Lemma \ref{lem:derivation}, we may assume that $O_{p'}(E)=1$.
So assume now that $p$ divides $|F(E)|$, and 
let $V$ be a minimal normal subgroup of $G$ contained in $Z(O_p(E))$. Our inductive hypothesis  applied to $G/V$ (or our assumption on $S$ if $E=V$) shows that there is a $p$-element $x$ of $G$ with $x\pi=z$ and $[D_0,x]\le V$. Then $d_0\rightarrow [x,d_0^{-1}]$ defines a derivation from $D_0$ to $V$. Since $D_0$ is a $p'$-group, this must be an inner derivation by the Schur-Zassenhaus theorem. That is, there exists $v\in V$ such that $[x,d_0^{-1}]=[v,d_0^{-1}]$ for all $d_0\in D_0$. It follows that $[D_0, v^{-1}x]=1$, whence the result.

Finally, assume that (3) holds, so that, in particular, $p$ divides $|E|$. In this case, we will prove, by induction on $|E|$, that some element of $E$ of order $p$ has centraliser projecting, modulo $E$, onto a subgroup of $S$ containing some $D\in\mathcal{S}$. Let $V$ be a minimal normal subgroup of $G$ contained in $Z(O_p(E))$. By Lemma \ref{lem:Z1}, we may assume that $V$ is a non-trivial $\mathbb{F}_p[S]$-module. If $V$ has a non-zero element with an insoluble centraliser, then the result clearly follows, so assume that $\mathrm{H}^2(S,V)=\mathrm{H}^1(D,V\hspace{-0.2cm}\downarrow_D)=0$ for all $D\in\mathcal{S}$. Since $E=\Phi(G)$, the condition $\mathrm{H}^2(S,V)=0$ implies that $O_p(E)\neq V$ (otherwise, $G/O_{p'}(E)$ would yield a non-split extension of $V$ by $S$). This completes the base case for induction. Applying our inductive hypothesis to $G/V$, we get that there is a non $G$-central $p$-element $x$ of $E$ and a subgroup $D\in\mathcal{S}$ with $[D,x]\le V$. Since $\mathrm{H}^1(D,V\hspace{-0.2cm}\downarrow_D)=0$, we deduce that the map $d\rightarrow [x,d^{-1}]$ is an inner derivation $D\rightarrow V$. The result now follows as in the previous paragraph.    
\end{proof}

\begin{corollary}\label{cor:RuleOut}
Let $G$ be a finite group in which the centralisers of all non-central elements are soluble. Then $G$ has no section isomorphic to $A_5\times A_6$, $A_5\times \PSp_4(3)$, $A_8$, $\Lg^{\pm}_5(2)$, $\Lg^{\pm}_5(3)$, 
$\PSp_6(2)$, $\PSp_6(3)$, 
$^{3}\mathrm{D}_4(2)$, or $^{3}\mathrm{D}_4(3)$.   
\end{corollary}
\begin{proof}
Suppose that $G$ has a section $X/Y\cong S$, where $S$ is one of the groups listed above. Choose $L$ minimal such that $LY=G$. Then $E:=L\cap Y\le \Phi(L)$, and $L/E\cong S$. Thus, $E=\Phi(L)$. Note also that each of the groups $S$ has the property that for some prime $p$, an element of $S$ of order $p$ has an insoluble centraliser. For example, if $S=A_8$, then we can take $p=3$; if $S=\Lg_5(2)$, we can take $p\in\{2,3\}$, and so on. Let $\mathcal{P}=\mathcal{P}(S)$ be the set of such primes $p$. 

Let $\mathcal{S}$ be the set of insoluble subgroups of $S$. We can use \textsc{Magma} \cite{Magma} to compute the set of irreducible $\mathbb{F}_p[S]$-modules, for each $p\in\mathcal{P}$. For each choice of $S$, we can check that for some prime $p\in\mathcal{P}(S)$, the hypothesis of Corollary \ref{cor:strategy} holds. This contradiction completes the proof.
\end{proof}

The three main consequences of this series of results are as follows. Recall the definitions of the sets $X_{\underline{\hspace{0.2cm}}}(\pi)$ from page 3.
\begin{proposition}\label{prop:Alt}
Let $G$ be a group of shape $G=E.S$ where $S=A_n$ with $n\geq 5$, and $E=\Phi(G)$. Let $\pi$ be a set of primes, and assume that all centralisers of non-central $\pi$-elements of $G$ are soluble. Then $n\le p+4+2\delta_{2,p}$ for all primes $p$ in $\pi$, and so $S \in X_{\mathrm{Alt}}(\pi)$.
\end{proposition}
 \begin{proof}
Suppose that the proposition fails, and let $G$ be a counterexample of minimal order. 
Then there eixts a prime $p\in\pi$ such that $n\geq m:=m(p)=p+5+2\delta_{2,p}$. Thus, the group $S$ contains a subgroup $L\cong A_m$. Let $\rho:G\rightarrow S$ be the natural projection, and let $X$ be a subgroup of $\rho^{-1}(L)$ which is minimal with the property that $EX=\rho^{-1}(L)$. Then $X=E_0.L$ with $E_0=\Phi(X)$. 
Thus, $G=X$ by the minimality of $G$ as a counterexample. Throughout, let $x$ be a $p$-cycle in $A_m$, and let $A$ be a subgroup of $C:=C_{A_m}(x)$ isomorphic to $A_5$. 

If $p>5$, then $m=p+5$, so $C$ contains the insoluble $p'$-subgroup $A$. The result therefore follows from Corollary \ref{cor:strategy}. Thus, we may assume that $p\le 5$, whence $S\in\{A_9\text{ }(p=2),A_8\text{ }(p=3),A_{10}\text{ }(p=5)\}$. We can then use \textsc{Magma} \cite{Magma} to construct all irreducible $\mathbb{F}_p[S]$-modules $V$. If $p=2$ or $p=5$, so that $S=A_9$ or $A_{10}$, then there are $9$ or $19$ possibilities for $V$ respectively, and in each case we can check that a non-zero vector has an insoluble stabiliser in $S$.

So we may assume that $p=3$, so that $S=A_{8}$.  
Let $D_1\cong A_5\le S$ be the pointwise stabiliser of a $3$-set; and let $D_2,D_3$ and $D_4$ be representatives of the three conjugacy classes of subgroups of $S\cong A_{8}$ isomorphic to $\Lg_3(2)$. We can check using \textsc{Magma} \cite{Magma} that some $D_i$, for $1\le i\le 4$, fixes a non-zero vector in $V$ unless $V$ is the unique irreducible $\mathbb{F}_3[S]$-module of dimension $90$. In this case, $\mathrm{H}^2(S,V)=H^1(D_i,V\downarrow_{D_i})=0$ for each $1\le i\le 4$. The result now follows from Corollary \ref{cor:strategy}. 
\end{proof}

\begin{proposition}\label{prop:PSp}
Let $G$ be a group of shape $G=E.S$ where $S=\PSp_n(q)$ with $n\geq 4$ even, and $E=\Phi(G)$. Let $\pi$ be a set of primes, and assume that all centralisers of non-central $\pi$-elements of $G$ are soluble. Then $n\le p+6+8\delta_{2,p}+6\delta_{3,p}+4\delta_{5,p}$ for all primes $p$ in $\pi$, and so $S \in X_{\mathrm{Class}}(\pi)$.
\end{proposition}
 \begin{proof}
Suppose that the proposition fails, and let $G$ be a counterexample of minimal order. Set $d:=d(p)=8+8\delta_{2,p}+6\delta_{3,p}+4\delta_{5,p}$.
Then there exists a prime $p\in\pi$ such that $n\geq p-1+d(p)$. 

The stabiliser in $S$ of a non-degenerate $(n-d(p))$-space has a Levi subgroup isomorphic to a central quotient of $\mathrm{GL}_{d(p)/2}(q)\times \mathrm{Sp}_{n-d(p)}(q)$.
Note that $p$ divides $q^{p-1}-1$ by Fermat's Little Theorem, and $n-d(p)\geq p-1$, so $\PSp_{n-d(p)}(q)$ contains a non-central element $x$ of order $p$. Also, $\mathrm{GL}_{d(p)/2}(q)$ contains a subgroup isomorphic to $A_9$, $A_8$, $\PSL_3(2)$ or $A_5$ for $p=2$, $p=3$, $p=5$ or $p>5$, respectively. Thus, $C:=C_S(x)$ contains such a subgroup, which we will call $A$.

If $p\geq 5$, then $(|A|,p)=1$, so $G$ contains a non-central element of order $p$ with an insoluble centraliser by Corollary \ref{cor:strategy}. If $p=2$ or $p=3$, we argue as in the proof of Proposition \ref{prop:Alt} 
\end{proof}

\begin{proposition}\label{prop:Spor}
Let $\pi$ be a set of primes, and let $G$ be a group of shape $G=E.S$ where $S$ is a sporadic simple $\pi$-group and $E=\Phi(G)$. Assume that all centralisers of non-central $\pi$-elements of $G$ are soluble. Then $S \in X_{\mathrm{Spor}}(\pi)$.
\end{proposition}
\begin{proof}
Note first that by Proposition \ref{prop:Alt}, the group $S$ cannot have a section isomorphic to an alternating group of degree larger than $p+4+2\delta_{2,p}$. 

Suppose first that $p=2$. Then we need to prove that $S\in\{\mathrm{M}_{11},\mathrm{M}_{22}\}$. As mentioned above, $S$ has no alternating section of degree larger than $8$. By Lemma \ref{lem:sporalt}, we deduce that $S\in\{\mathrm{M}_{11},\mathrm{M}_{12},\mathrm{M}_{22},\mathrm{M}_{23},\mathrm{M}_{24},\mathrm{HS},\mathrm{J}_2,\mathrm{Co}_2,\mathrm{Co}_3,\mathrm{McL},\mathrm{Suz},\mathrm{He},\mathrm{J}_1,\mathrm{O'N},\mathrm{J}_3,\mathrm{Ru},\mathrm{J}_4\}$.
By Corollary \ref{cor:strategy}, we may assume that $|E|$ is even, and that there is an irreducible $\mathbb{F}_2[S]$-module in which all non-zero elements have soluble centralisers. One can use \textsc{Magma} \cite{Magma} to check that no such module exists if $S\in\{\mathrm{M}_{12},\mathrm{M}_{23},\mathrm{M}_{24},\mathrm{HS},\mathrm{J}_2,\mathrm{McL},\mathrm{He},\mathrm{J}_1,\mathrm{J}_3,\mathrm{Ru}\}$. Since $\mathrm{Co}_2,\mathrm{Co}_3,\mathrm{Suz},\mathrm{O'N}$ and $\mathrm{J}_4$ all contain a section isomorphic to one of $\mathrm{M}_{12},\mathrm{M}_{23},\mathrm{M}_{24}$, $\mathrm{HS},\mathrm{J}_2,\mathrm{McL},$ or $\mathrm{J}_1$, the result follows.

Suppose next that $p=3$. We need to show that $S\in\{\mathrm{M}_{11},\mathrm{M}_{12},\mathrm{M}_{22},\mathrm{J}_1\}$. So assume that $S$ is not one of these groups. We can also assume that $S$ has no alternating section of degree bigger than $7$. The possibilities for $S$ are then given by Lemma \ref{lem:sporalt}. Thus, $S$ is one of
$$\mathrm{J}_{2},\mathrm{Suz},\mathrm{He},\mathrm{O'N},\text{ or }\mathrm{J}_{3}.$$
As above, we may assume that $|E|$ is divisible by $3$, and that there is an irreducible $\mathbb{F}_3[S]$-module in which all non-zero elements have soluble centralisers.
One can use \textsc{Magma} \cite{Magma} to check that no such module exists when $S$ is $\mathrm{J}_{2},\mathrm{He},\mathrm{O'N},\text{ or }\mathrm{J}_{3}$. Since $\mathrm{Suz}$ contains a section isomorphic to $J_2$, the result follows.

The cases $p=5,7$ are argued in exactly the same way, so lets move on to the cases $p\geq 11$. An easier argument works here. We exhibit the case $p=11$: in this case, the only group we need to rule out is the monster group $\mathrm{M}$. So assume that $S=\mathrm{M}$. Then $S$ has an element of order $11$ with centraliser containing the $11'$-group $A_6$. The result then follows from Corollary \ref{cor:strategy}. The same approach deals with all primes $p>11$.   
\end{proof}

We can now prove Theorem \ref{thm:Main}.
\begin{proof}[Proof of Theorem \ref{thm:Main}] 
Let $G$ be a counterexample to the theorem of minimal order, with $\pi$ fixed. 

Let $N$ be a maximal normal subgroup of $G$. If $M$ is a maximal subgroup of $G$ not containing $N$, then minimality implies that all non-abelian composition factors of $M$ and $N$ lie in the required list. Since $MN=G$, the result follows.

Thus, we may assume that $G=E.S$ with $E=\Phi(G)$ and $S$ a non-abelian simple group. The choice of $G$ implies that
$$ S \notin X_{\mathrm{Alt}}(\pi)\cup X_{\mathrm{Class}}(\pi)\cup X_{\mathrm{Exc}} \cup X_{\mathrm{Spor}}(\pi).$$ In particular $S$ is not exceptional.
 If $S$ is an alternating or sporadic group, then from Propositions \ref{prop:Alt} and \ref{prop:Spor} we deduce $S \in X_{\mathrm{Alt}}(\pi)\cup X_{\mathrm{Spor}}(\pi)$, a contradiction.

So assume that $S$ is a classical group of Lie type. Let $p\in\pi$, and suppose that $S$ has a subgroup $S_0$ with $S_0\cong A_{p+5+2\delta_{2,p}}$. Choose $G_0\le G$ minimal with the property that $G_0$ projects onto $S_0$ modulo $E$. Then $E_0:=G_0\cap E\le \Phi(G_0)$, and $G_0$ is a group of shape $E_0.S_0$. The required contradiction then follows from Proposition \ref{prop:Alt}.    

Thus, $S$ has no alternating subgroup of degree greater than $p+4+2\delta_{2,p}$.
 Hence, $S \in X_{\mathrm{Class}}(\pi)$ by Lemma \ref{lem:AltEmbeddings}

Finally, the statement on simple groups follows from Lemma \ref{lem:simple}.
\end{proof}

\section{Proofs of Theorem \ref{thm:MainClass}, and Corollaries \ref{cor:MainClass} and  \ref{cor:NCgraph}}

We begin by proving that finite groups with $R(G)=1$ often have a proper insoluble centraliser.
\begin{lemma}\label{lem:KeyFactors}
Let $G$ be a finite group with $R(G)=1$, and let $N=S^t$ be a minimal normal subgroup of $G$, with $S$ a non-abelian simple group.
Suppose that $|G|$ has a prime divisor $p$ with $p\nmid |N_G(S)/C_G(S)|$. Then $G$ has an element $x$ of order $p$ and an insoluble subgroup $D$ of $p'$-order such that $[D,x]=1$.
\end{lemma}
\begin{proof}
Let $x$ be an element of $G$ of order $p$. Clearly, we may assume that $x$ does not centralise any of the direct factors in $N$. Then since $p\nmid|N_G(S)/C_S(S)|$, we must have that $p$ divides $t$, and there is an embedding $N\langle x\rangle\hookrightarrow N:\langle \sigma\rangle$, where $\sigma\in S_{t}$ acts as a fixed point free permutation of order $p$ on the direct factors in $N$. In particular, $x$ centralises some diagonal subgroup of $N$. Taking $D$ to be this diagonal subgroup then yields the result.
\end{proof}

We immediately deduce the following corollary, which shows that the structure of an insoluble group in which all proper centralisers are soluble is very restricted.
\begin{corollary}\label{cor:KeyFactors}
Let $G$ be a finite insoluble group in which the centraliser of every non-central element is soluble, and let $N/R(G)=S^t$ be a minimal normal subgroup of $G/R(G)$, where $S$ is a non-abelian simple group. Then $G/R(G)$ is a $\pi$-group, where $\pi=\pi(\Aut(S))$. 
\end{corollary}

Recall that for a $p$-group $A$, $\Omega_1(A):=\langle x\text{ : }|x|=p\rangle$.
\begin{lemma}\label{lem:centr}
 Let $p$ be prime, let $A$ be a finite abelian $p$-group, and let $\alpha$ be a $p'$-automorphism of $A$. If $[\Omega_1(A),\alpha]=1$, then $[A,\alpha]=1$. In particular, if $X$ is an insoluble subgroup of $\Aut(A)$ and $[\Omega_1(A),X]=1$, then $[A,X]=1$.
\end{lemma}
\begin{proof}
Suppose that $[\Omega_1(A),\alpha]=1$, and write $p^k$ for the exponent of $A$. We will prove, by induction on $k$, that $A$ is centralised by $\alpha$. The case $k=1$ follows by hypothesis. So assume $k>1$, and let $x\in A$. Then $1=[\alpha,x^p]$ by the inductive hypothesis, so $1=[\alpha,x]^p$, since $A$ is abelian and $\alpha$-invariant. Thus, $[\alpha,x]\in \Omega_1(A)$, whence $[x,\alpha,\alpha]=1$. Since $(|\alpha|,p)=1$, we deduce that $[x,\alpha]=1$. The claim follows.

 The last part of the lemma follows since if $X$ is insoluble, then $\langle \alpha^X\rangle$ is insoluble for some $p'$-element $\alpha$ of $X$. 
\end{proof}

In our next lemma, for a set of primes $\pi$, we write $R_{\pi}(G)$ for the largest normal $\pi$-soluble subgroup of $G$; and $O_{\pi}(G)$ for the largest normal $\pi$-subgroup of $G$. Recall also that a \emph{Frattini central cover} of a group $G$ is a group $X$ with a normal subgroup $N\le \Phi(X)\cap Z(X)$ satisfying $X/N\cong G$.
\begin{lemma}\label{lem:pigroupslift}
Let $\pi$ be a set of primes and let $G$ be a finite group of shape $G=R_{\pi}(G).A$ where $A$ is almost simple. Let $\pi_1:=\pi(A)\cap \pi$. 
Suppose that
\begin{enumerate}
    \item[\upshape(i)] every Sylow $r$-subgroup of $G$, for $r\in \pi_1$, is either $r$-central or has nilpotency class equal to the class of a Sylow $r$-subgroup of $A$; and
    \item[\upshape(ii)] every non-central $\pi_1$-element of $G$ has a soluble centraliser. 
\end{enumerate}
Then $O_{\pi_1'}(G)$ is soluble, and $G/Z(G)=D.A$ for some soluble $\pi_1'$-group $D$.   
\end{lemma}
\begin{proof}
We will prove the claim by induction on $|R|$, where $R=R_{\pi}(G)$. Note that 
$\soc(A)$ is not a $\pi_1'$-group, since $G/R_{\pi}(G)$ is almost simple and $G\neq R_{\pi}(G)$.
Set $D_0:=O_{\pi_1'}(G)$. Then $D_0$ is contained in $R$, since $G/R$ is almost simple and $\soc(A)$ is not a $\pi_1'$-group. 

If $R=1$, then the result is trivial, and this can serve as the base case for induction. So assume that $R>1$. We will first prove that $D_0$ is soluble. Indeed, suppose otherwise, and let $N$ be a minimal normal subgroup of $G$ contained in $D_0$. If $N$ is soluble, then since $D_0$ is a $\pi_1'$-group, our inductive hypothesis can be applied to $G/N$. Thus, we may assume that $N=S^t$ with $S$ a non-abelian simple $\pi_1'$-group. In a slight abuse of notation, we will identify $S$ with one of the direct factors in $N$. Let $p\in\pi_1$ with $p$ dividing $|\soc(A)|$, and write $\Omega_1(G)$ for the subgroup of $G$ generated by the elements of $G$ of order $p$. Then $C_G(S)/Z(G)$ is a $\pi_1'$-group, since the centraliser of every non-central $\pi_1$-element in $G$ is soluble. If $\Omega_1(G)\le N_G(S)$, then it would follow that $\Omega_1(G)\le C_G(S)$, since $N_G(S)/C_G(S)$ is an almost simple group with $p'$-socle. But then $\Omega_1(G)\le Z(G)$, which is a contradiction, since $p$ divides $|\soc(A)|$. Thus, there exists an element $x$ of $G$ of order $p$ which acts fixed point freely on the $t$ factors in $N$. In particular, $G$ has a subgroup $X$ of the form $X:=S^p\langle x\rangle$. Thus, by Lemma \ref{lem:KeyFactors}, $X$ has a $p$-element $g$ with $C_X(g)$ containing an insoluble subgroup of $p'$-order. This is a contradiction, by Corollary \ref{cor:strategy} (applied with Hypothesis (2)). Thus, $D_0$ is soluble, as needed. 

Next, we claim that
\begin{align}\label{lab:claim:GDA}
 G/D_0=(Z(G)D_0/D_0)\circ \hat{A}   
\end{align}
with $Z(G)D_0/D_0$ a $\pi_1$-group, where $\hat{A}$ is a Frattini central cover of the almost simple group $A$ in the statement of the lemma. 
This will suffice to prove the lemma. Indeed, since $(|D_0|,|Z(G)D_0/D_0|)=1$, we have $Z(G)D_0/D_0=F/D_0$ with $F=Z\times D_0$ and $Z\le Z(G)$. Writing $B/D_0=\hat{A}$, we then have $BZ=G$, so $G/Z=B/B\cap Z\cong D.A$ where $D$ is a soluble $\pi_1'$-group.

To prove \eqref{lab:claim:GDA},  write bars to denote reduction modulo $D_0=O_{\pi_1'}(G)$. Since the generalised Fitting subgroup of a finite group is self-centralising, it will suffice to prove that $F^*(\ol{G})=\ol{Z(G)}\circ \hat{S}$ where $\hat{S}$ is a Frattini central cover of a non-abelian simple group $S$. Note that $R_{\pi}(\ol{G})=R/D_0=R_{\pi}(G)/D_0$ since $D_0$ is a $\pi'$-group. Thus, $\ol{G}/R_{\pi}(\ol{G})\cong A$. 

To prove that $F^*(\ol{G})=\ol{Z(G)}\circ \hat{S}$, we will begin by showing that  
\begin{align}\label{lab:k}
    F(\ol{G})\le \ol{Z(G)}.
\end{align}
To see this, let $r$ be a prime with $O_r(\ol{G})\neq 1$. Then $r\in \pi_1$, since $O_{\pi_1'}(\ol{G})=1$. We need to prove that
$O_r(\ol{G})\le \ol{Z(G)}$. Suppose first that the nilpotency class of a Sylow $r$-subgroup of $A$ equals the nilpotency class, say $c$, of a Sylow $r$-subgroup of $G$.  Then the same is true in $\ol{G}$, so $O_r(\ol{G})$ is centralised by $\gamma_c(Q)D/D$ for all Sylow $r$-subgroups $Q$ of $G$. Thus, $C_{\ol{G}}(O_r(\ol{G})R_{\pi}(\ol{G})/R_{\pi}(\ol{G})$ is non-trivial, so $C_{\ol{G}}(O_r(\ol{G})$ is insoluble, since $\ol{G}/R_{\pi}(\ol{G})$ is almost simple. It follows from Hypothesis (ii) of the lemma, and Lemma \ref{lem:derivation}, that $O_{r}(\ol{G})$ is generated by a set of elements of $\ol{Z(G)}$, i.e. $O_{r}(\ol{G})\le \ol{Z(G)}$.

So we may assume that a Sylow $r$-subgroup of $G$ is $r$-central. Since $D$ is $r'$-group, the same is true in $\ol{G}$. Let $E$ be an $r$-subgroup of $G$ with $\ol{E}$ a characteristic abelian $r$-subgroup of $F(\ol{G})$. Then $\Omega_1(\ol{E})$ is an elementary abelian subgroup of $\Omega_1(O_r(\ol{G}))$. Let $X$ be the subgroup of $G$ generated by all of the Sylow $r$-subgroups of $G$. Then $XR/R\geq \soc(G/R)$, since $X\unlhd G$ and $G/R$ is almost simple with socle of order divisible by $r$. Also, $X\le C_G(\Omega_1(\ol{E}))$ since $\ol{G}$ is $r$-central. It follows from Lemma \ref{lem:centr} that $[\ol{E},X]=1$. Thus, arguing as in the paragraph above, we see that for each $e\in E$, there exists a subgroup $X_0$ of $G$ with $X_0D/D=XD/D$, and $[e,X_0]=1$. Since the centralisers of all non-central $r$-elements of $G$ are soluble, we deduce that $E\le Z(G)$. Thus, all characteristic abelian $r$-subgroups of $F(\ol{G})$ are in $\ol{Z(G)}$. It follows that $O_r(\ol{G})$ is either abelian or extraspecial. But extraspecial groups are not $r$-central. Thus, $O_r(\ol{G})$ is abelian, whence contained in $\ol{Z(G)}$. This completes the proof of \eqref{lab:k}. 

Since $\ol{G}\neq F(\ol{G})$ and $C_{\ol{G}}(F^*(\ol{G}))\le F^*(\ol{G})$, we deduce that $F(\ol{G})\neq F^*(\ol{G})$, whence $\ol{G}$ has a normal subgroup $T$ which is a central product of (say) $t$ isomorphic copies of a quasisimple group $S$. If $t>1$, then since $O_{\pi_1'}(\ol{G})=1$, there exists a prime $r\in\pi_1$, and a non-central element in $T$ of order $r$ with an insoluble centraliser. Applying Lemma \ref{lem:derivation} then yields a non-central $r$-element of $G$ with an insoluble centraliser, contrary to assumption. So we must have $t=1$. The claim \eqref{lab:claim:GDA} follows, whence the result.          
\end{proof}

\begin{proof}[Proof of Theorem \ref{thm:MainClass}]
We may assume throughout that $G$ is not $\pi$-soluble. We will first prove that 
\begin{align}\label{lab:pipclaim}
    G\text{ has shape }G=R_{\pi}(G).A
\end{align}
where $A$ is an almost simple group in which every non-trivial $\pi$-element has soluble centraliser. 

Let $G$ be a counterexample to \eqref{lab:pipclaim} of minimal order. Suppose first that $G$ has trivial $\pi$-soluble radical. Let $N$ be a minimal normal subgroup of $G$, so that $N=S^t$ for $S$ a non-abelian simple group with at least one prime divisor in $\pi$. Since all non-central $\pi$-elements have soluble centralisers, we must have $t=1$ and $C_G(S)=1$. Thus, $G$ is almost simple. The claim follows.

Next, assume that $G$ has non-trivial $\pi$-soluble radical $R:=R_{\pi}(G)$. Suppose that there exists a maximal subgroup $M$ of $G$ with $R\not\le M$. Then $G=MR$, and $M\cap R\le R_{\pi}(M)$. By the minimality of $G$ as a counterexample, we have that either $M$ is $\pi$-soluble, or $M/R_{\pi}(M)$ is almost simple in which all non-trivial $\pi$-elements have soluble centraliser. Since $G=MR$, we deduce that either $G$ is $\pi$-soluble, or $M\cap R=R_{\pi}(M)$ and $G/R\cong M/M\cap R$. This gives us what we need.  

Thus, we may assume that $R\le \Phi(G)$, whence $R=\Phi(G)$. In particular, $R$ is nilpotent. Let $E/R$ be the socle of $G/R$, so that $E/R$ is a direct product of non-abelian simple groups, each of which is divisible by at least one prime in $\pi$. We will first prove that 
\begin{align}\label{claim:ER}
    \text{$E/R$ is simple.}
\end{align}
To do this, assume that $E/R$ is not simple, and fix $r\in\pi\cap \pi(E/R)$. Then there is an $r$-element $x$ of $E/R$ with an insoluble centraliser. If $r$ does not divide $|R|$, then Lemma \ref{lem:derivation} shows that there is a non-central $r$-element of $G$ with an insoluble centraliser, contrary to hypothesis. 

So we may assume that $r$ divides $|R|$. Let $V$ be a minimal subgroup of $G$ contained in $\Omega_1(Z(O_r(R)))\le Z(R)$. We claim that $V\le Z(G)$. Indeed, recall that $G$ satisfies property $(\ast)_{\pi}$. Suppose first that a Sylow $r$-subgroup of $G$ is $r$-central. Then all $r$-elements of $G$ centralise $V$. Since $E/R$ has a non-abelian simple direct factor of order divisible by $r$, it follows that $V$ has an insoluble centraliser. Since all non-central $r$-elements of $G$ have soluble centraliser, we deduce that $V\le Z(G)$. 

Suppose next that $G$ has a non-abelian simple section $S=L/K$ with a Sylow $r$-subgroup $P/K\le S$ of nilpotency class $c$, where $c$ is the nilpotency class of a Sylow $r$-subgroup of $G$. Then for each $\ell\in L$, the group $\gamma_c(P^{\ell})$ is not contained in $K$, and centralises $O_r(R)$. Thus, $S$ centralises $O_r(R)\geq V$. As above, our hypothesis then implies that $V\le Z(G)$.    

So we have shown that $V\le Z(G)$. It follows from Lemma \ref{lem:Z1} that $E/V$ has no $r$-element with an insoluble centraliser. We can then apply the argument above to $E/V$ again. Doing this repeatedly in fact yields that all chief factors of $G$ contained in $O_r(R)$ are central. The required contradiction then follows from Corollary \ref{cor:Z1}. 
Thus, \eqref{claim:ER} holds.

We therefore have that $G/R_{\pi}(G)$ is almost simple. It now follows from Lemma \ref{lem:pigroupslift} that $O_{\pi_1'}(G)$ is soluble, and $G/O_{\pi_1'}(G)$ has the stated shape.   

Thus, all that remains is to prove that every non-trivial $\pi$-element of $A=G/O_{\pi_1'}(G)$ has soluble centraliser.
To do this, let $r\in\pi$ be a prime divisor of $|A|$, and suppose that $S$ has an $r$-element $x$ with an insoluble centraliser. If $r$ does not divide $|R|$, then we get a contradiction from Lemma \ref{lem:derivation}. So assume that $r$ divides $|R|$. Then arguing as above, we see that $Z(G)$ is non-trivial. We then get a contradiction from Lemma \ref{lem:Z1}. It follows that if $r\in\pi$, then all $r$-elements of $A$ have soluble centralisers.
\end{proof}

\begin{proof}[Proof of Corollary \ref{cor:MainClass}]
By Theorem \ref{thm:MainClass}, we see that $G=R(G).A$, where $A$ is almost simple and every involution in $A$ has a soluble centraliser. The corollary now follows immediately from Lemma \ref{lem:simple}.     
\end{proof}

\begin{proof}[Proof of Corollary \ref{cor:NCgraph}] Let $G$ and $H$ be as in the statement of the theorem, and let $\varphi:\Gamma_G\rightarrow \Gamma_H$ be an isomorphism. Note that $|G|=|H|$ by \cite[Proposition 3.14]{AAM}.

For $x\in G\setminus Z(G)$, let $\Lambda(G,x)$ be the subgraph of $G$ induced by the vertices not adjacent to $x$. Then let $\Lambda^*(G,x)$ be the subgraph of $\Lambda(G,x)$ with the isolated vertices removed. Note that $\Lambda^*(G,x)$ is isomorphic to the non-commuting graph of $C_G(x)$. Further, $\varphi(\Lambda^*(G,x))$ is isomorphic to the non-commuting graph of $C_H(\varphi(x))$. It follows from the minimality of $G$ that $C_G(x)$ is soluble. Part (i) now follows from Theorem \ref{thm:2ndMain}. 

Finally, by Flavell's theorem \cite{Flavell}, $G$ has a $2$-generated insoluble subgroup $T=\langle x,y\rangle$. If $C_G(x)\cap C_G(y) \neq Z(G)$, then $G$ would have a non-central element with an insoluble centraliser, contrary to the above. Thus, $C_G(x)\cap C_G(y)=Z(G)$. It follows that $a:=\varphi(x)$ and $b:=\varphi(y)$ have the same property in $H$. 
\end{proof}

\section{Proof of Theorem \ref{thm:2ndMain}}

The purpose of this section is to prove Theorem \ref{thm:2ndMain}. Recall that
$$\mathcal{Q}:=\{2^{p_1},3^{p_2},p\mid p_i,p\text{ odd primes, }p\equiv 0,\pm 2\pmod{5}\}.$$
We begin with a technical lemma.

\begin{lemma}\label{lem:simpleprimes}
Let $S_1,\hdots,S_m$ be finite simple groups with $m\geq 2$. Suppose that for all $i\neq j$ we have $\pi(S_i)=\pi(D_0)$ and $\pi(S_j)=\pi(E_0)$ for all insoluble sections $D_0$ of $S_j$, $E_0$ of $S_i$. Then one of the following holds.
\begin{enumerate}
    \item[\upshape(I)] $S_1=\hdots=S_m={}^2\mathrm{B}_2(2^p)$ for some odd prime $p$.
    \item[\upshape(II)] $S_i\in\{\Lg_2(q)\text{ : }q\in\mathcal{Q}\}$ for all $i$. Furthermore, if $S_i=\Lg_2(q_1)$ and $S_j=\Lg_2(q_2)$ for some $q_1\neq q_2$, then $(q_1,q_2)=1$.  
    \item[\upshape(III)] $S_i\in\{A_5,A_6,\PSp_4(3)\}$ for all $i$;
    \item[\upshape(IV)] $S_i\in\{\Lg_3(2),\Lg_2(8),\Ug_3(3)\}$ for all $i$; or
    \item[\upshape(V)] $S_i=\Lg_3(3)$ for all $i$.
\end{enumerate}
\end{lemma}
\begin{proof}
It suffices to prove the lemma with $m=2$. For ease of notation write $S_1=S$ and $S_2=T$.
Our assumptions of course imply $\pi(S)=\pi(T)$. Suppose first that either $S$ or $T$ contains a non-abelian simple subgroup $X$ with at most three distinct prime divisors. Then our assumptions force $|\pi(S)|=|\pi(T)|=3$. Further, by \cite[Lemma 2.1]{ZCCL}, $S,T\in\{A_5, A_6, \Lg_2(8), \Lg_3(2)$, $\Lg_3(3), \PSp_4(3), \Ug_3(3), \Lg_2(17)\}$.
These are cases (III)--(V) in the statement of the lemma.

Thus, we may assume that neither $S$ nor $T$ contains a non-abelian simple subgroup with at most three distinct prime divisors. It is well-known that all sporadic groups have a subgroup isomorphic to $A_5$. It follows, in particular, that neither $S$ nor $T$ is an alternating or sporadic group. So for the remainder of the proof, we may assume that $S$ and $T$ are groups of Lie type.

If $S$ [or $T$] is classical, then it is shown by Holt in \cite{MathOverflow} that $S$ [or $T$] contains a subgroup isomorphic to $A_5$ unless possibly $S$ [or $T$] lies in $\{\Lg_2(q),\Lg_3(q),\Ug_3(q)\}$ and $q(q-1)(q+1)$ is not divisible by $60$.

Thus, we may assume that $S$ and $T$ lie amongst the exceptional groups of Lie type together with the groups $\{\Lg_2(q),\Lg_3(q),\Ug_3(q)\}$ with $q(q-1)(q+1)$ not divisible by $60$ ($\ast$).
If $S$ and $T$ are both minimal simple groups, then Thompson's classification \cite[Corollary 1]{Thompson}, together with Zsigmondy's theorem, show that either (I), (II) or (V) must hold.

So we may assume that $S$ is not a minimal simple group. Then $S$ contains a minimal simple group $S_0$ by \cite{BarryWard}. Let $M$ be a maximal subgroup of $S$ containing $S_0$. Thus, $\pi(M)\subseteq \pi(\Aut(T))$ and $\pi(T)\subseteq \pi(\Aut(S_0))$ by assumption. Then $\pi(M)=\pi(S)$ by assumption, so the possibilities for $(S,M)$ are given by \cite[Table 10.7]{LPS}. The only possibilities satisfying ($\ast$) and not minimal simple are $S\in\{\mathrm{G}_2(3),{}^2\mathrm{F}_4(2)'\}$. Both of these groups have insoluble subgroups $D$ with $\pi(D)\neq \pi(S)$. This completes the proof. 
\end{proof}

\begin{corollary}\label{cor:simpleprimes}
Let $G$ be a finite simple group, and suppose that each of the following holds:
\begin{enumerate}
    \item[\upshape(i)] $G$ has no alternating sections of degree greater than $7$; and
    \item[\upshape(ii)] all sections $X$ of $G$ satisfy the following: For all non-central elements $x$ of $X$ of prime order $p$, the centraliser $C_X(x)$ has no insoluble subgroups of $p'$-order.
\end{enumerate}
Then setting $\mathcal{Q}:=\{2^{p_1},3^{p_2},p\mid p_i,p\text{ odd primes, }p\equiv 0,\pm 2\pmod{5}\}$, $G$ is one of the following.
\begin{enumerate}
    \item[\upshape(1)] $\mathrm{M}_{11}$, $\mathrm{M}_{22}$, $A_5$, $A_6$, $A_7$, $\mathrm{G}_2(3)$, $\Lg^{\pm}_5(2)$, $\Lg^{\pm}_5(3)$, $\Lg^{\pm}_6(2)$, $\Lg^{\pm}_6(3)$, $\PSp_4(2)$, $\PSp_4(3)$, $\PSp_6(2)$, $\PSp_6(3)$, $\PSp_8(2)$, $\PSp_8(3)$, ${}^2\mathrm{F}_4(2)'$, $^{3}\mathrm{D}_4(2)$, $^{3}\mathrm{D}_4(3)$.
    \item[\upshape(2)] $\Lg^{\pm}_4(q)$, $\PSp_4(q)$, $\mathrm{G}_2(q)$ or ${}^2\mathrm{F}_4(q)'$ with $q\in\mathcal{Q}$.
    \item[\upshape(3)] $\Lg_2(q)$, $\Lg^{\pm}_3(q)$,  ${}^2\mathrm{B}_2(q)$ or ${}^2\mathrm{G}_2(q)$.
\end{enumerate}       
\end{corollary}
\begin{proof}
The case where $G$ is an alternating group is of course trivial by (i). If $G$ is sporadic, then Lemma \ref{lem:sporalt} implies that $G$ is one of $\mathrm{M}_{11}$, $\mathrm{M}_{12}$, $\mathrm{M}_{22}$, $\mathrm{J}_{2}$, $\mathrm{Suz}$, $\mathrm{He}$, $\mathrm{J}_{1}$, $\mathrm{O'N}$ or $\mathrm{J}_3$.   
The result can then be checked using the Atlas of finite groups \cite{ATLAS}.

So we may assume that $G$ is a finite group of Lie type. By Assumption (i) and Lemma \ref{lem:AltEmbeddings}, $G$ is either classical with natural module of dimension at most $7$; isomorphic to $\PSp_n(q)$ with $q$ odd and $n\le 14$; or $G$ is ${}^2\mathrm{B}_2(q)$, $\mathrm{G}_2(q)$,
${}^2\mathrm{G}_2(q)$, ${}^2\mathrm{F}_4(q)'$ or $^{3}\mathrm{D}_4(q)$. 

By Assumption (ii), any section of $G$ of the form $S\times T$, with $S$ and $T$ non-abelian simple groups, must satisfy the conclusion of Lemma \ref{lem:simple}. That is, $S$ and $T$ must be one of the cases (I)--(V) in Lemma \ref{lem:simpleprimes}.

These restrictions give us what we need in all cases. Indeed,
suppose first that $G=\Lg^{\epsilon}_n(q)$, for $\epsilon\in\{\pm\}$. By Lemma \ref{lem:AltEmbeddings}, $n\le 6$. Thus, we may assume that $q>3$. For all partitions $a+b=n$ of $n$, the group $G$ has a section $X$ isomorphic to $\Lg^{\epsilon}_a(q)\times\Lg^{\epsilon}_b(q)$. 
If $n\geq 5$, then taking $a=n-2$ and $b=2$ above gives us what we need by taking $p$ to be a Zsigmondy prime divisor of $q^{n-2}-(\epsilon)^{n-2}$, since $\Lg_2(q)\cong \Lg_2^-(q)$ is then a $p'$-group centralising a $p$-element of $\Lg^{\epsilon}_{n-2}(q)$. So assume further that $n\le 4$. If $n=4$, then the groups $S$ and $T$ are as in one of the cases (I)--(IV) in Lemma \ref{lem:simpleprimes}. It follows that $q\in\{2^{p_1},3^{p_2},p\text{ : }p_i\text{ an odd prime, }5\mid p^2+1\text{ }\}$. This completes the proof in this case. 

By Lemma \ref{lem:sporalt}, the only remaining possibility when $G$ is classical is $G=\PSp_n(q)$ with $n\le 14$ for $q$ odd; and $n\le 6$ for $q$ even. Then $G$ has a section $\PSp_a(q)\times \PSp_b(q)$ for $a,b$ even and $a+b=n$. We then argue as above. The exceptional cases are entirely similar. For example, ${}^3\mathrm{D}_4(q)$ has a section isomorphic to $\Lg_2(q^3)\times \Lg_2(q)$, so Lemma \ref{lem:simpleprimes} implies that $q$ is $2$ or $3$.
\end{proof}

We are finally ready to prove Theorem \ref{thm:2ndMain}.
\begin{proof}[Proof of Theorem \ref{thm:2ndMain}]
Let $G$ be a counterexample to the theorem of minimal order. Thus, every non-central element of $G$ has soluble centraliser, but $G$ is not soluble, and $G$ is not as in parts (ii) and (iii) of the theorem. Then by Lemma \ref{lem:centr}, $Z(G)=1$. Thus, all non-trivial elements of $G$ have soluble centraliser. For the remainder of the proof, let $R:=R(G)$ be the soluble radical of $G$.

Suppose first that $G/R$ is almost simple, with socle $S$.
Assume that $S$ has a section $X:=L/K$ such that, for some non-central element $x\in X$ of prime order $p$, $C_X(x)$ contains an insoluble subgroup of $p'$-order. Choose $Y\le L$ minimal with the property that $YK=L$. Then $E:=Y\cap K\le \Phi(Y)$, and $Y/E\cong X$. Corollary \ref{cor:strategy} (applied with Hypothesis (2)) then implies that $L$ has an element of order $p$ with an insoluble centraliser -- a contradiction. We deduce that no such section exists. In the same way, we can apply Proposition \ref{prop:Alt} to deduce that $S$ has no alternating sections of degree greater than $7$. It follows that all non-abelian composition factors of $G$ lie in Corollary \ref{cor:simpleprimes}(1)--(3). The groups in Part (1) of the latter list which do not lie in Theorem \ref{thm:2ndMain}(ii)(1) can all be ruled out using Corollary \ref{cor:RuleOut}.

Thus, $G/R(G)$ is not almost simple. We need to prove that for one of the sets $\pi$ in Table \ref{tab:Main}, $G/R(G)$ is a $\pi$-group and all non-abelian composition factors of $G$ lie in $S(\pi)$. 
To this end, note first that since $G/R$ is not almost simple, we have $\soc(G/R)=S_1\times\hdots\times S_m$ with $S_i$ non-abelian simple for all $i$, and $m\geq 2$. As in the first paragraph above, Corollary \ref{cor:strategy} then forces $\pi(S_i)=\pi(D_0)$ and $\pi(S_j)=\pi(E_0)$ for all insoluble sections $D_0$ of $S_j$, $E_0$ of $S_j$, and all $i,j$. Thus, the $S_i$ satisfy Lemma \ref{lem:simpleprimes}(I)--(V).

Suppose that case (II) of Lemma \ref{lem:simpleprimes} occurs. That is, each $S_i$ has the form $S_i=\Lg_2(q_i)$ for some prime power $q_i$, with either $(q_i,q_j)=1$ or $q_i=q_j$, for each $i,j$. Let $\mathcal{Q}_0:=\{q_1,\hdots,q_m\}\subseteq \mathcal{Q}$. Since $\pi(S_i)=\pi(S_j)$ for all $i,j$, it then follows from the definition of $\mathcal{Q}$ that $|\mathcal{Q}_0|\le 3$. 

Note also that if case (III) of Lemma \ref{lem:simpleprimes} occurs, then $\pi(S_i)=\{2,3,5\}$ for all $i$, and only one of $A_5,A_6$ or $\PSp_4(3)$ can occur by Corollary \ref{cor:RuleOut}.


Now, let $\pi:=\bigcup_{i=1}^m\pi(\Aut(S_i))$, and let $K/R$ be the kernel of the action of $G/R$ on the set $\{S_1,\hdots,S_m\}$. Then $\pi$ lies in Table 1, and $(K/R)/\soc(G/R)$ is a soluble $\pi$-group. Thus, all that remains is to prove that $G/K$ is a $\pi$-group, and that all non-abelian composition factors of $G/K$ lie in $S(\pi)$. 

So let $p$ be a prime divisor of $|G/K|$, and assume that $p\not\in\pi$. Let $Rx\in G/R$ of order $p$. Then $Rx\not\in K/R$, since $K/R$ is a $\pi$-group. Thus, $Rx$ has a non-trivial orbit on $\{S_1,\hdots,S_m\}$. It follows that $G/R$ has a section $X\cong S\wr C_p$.  It follows from Lemma \ref{lem:KeyFactors} that $X$ has an element $x$ of order $p$ such that $C_{X}(x)$ has an insoluble subgroup of $p'$-order. This contradicts the second paragraph above. The proof is complete.
\end{proof}

\section*{Acknowledgments}
The first two authors are members of the \emph{Gruppo Nazionale per le Strutture Algebriche, Geometriche e le loro Applicazioni (GNSAGA -- INdAM)}. 
The authors wish to express their gratitude to the \emph{Centro Internazionale per la Ricerca Matematica (CIRM)} in Trento, where this research was carried out. 
This work was partially supported by GNSAGA -- INdAM. 
Moreover, the second author was funded by the European Union -- Next Generation EU, Missione~4 Componente~1, CUP~B53D23009410006, PRIN~2022 (\emph{Group Theory and Applications}, code~2022PSTWLB).

\bibliographystyle{amsplain}
\bibliography{books.solublecentralizer}
\end{document}